\documentclass[11pt]{article}
\usepackage{amsfonts}
\usepackage{mathrsfs}
\usepackage{cite}
\usepackage{amssymb,amsmath} 
 \allowdisplaybreaks

\catcode`!=11
\let\!int\int \def\int{\displaystyle\!int}
\let\!lim\lim \def\lim{\displaystyle\!lim}
\let\!sum\sum \def\sum{\displaystyle\!sum}
\let\!sup\sup \def\sup{\displaystyle\!sup}
\let\!inf\inf \def\inf{\displaystyle\!inf}
\let\!cap\cap \def\cap{\displaystyle\!cap}
\let\!max\max \def\max{\displaystyle\!max}
\let\!min\min \def\min{\displaystyle\!min}
\let\!frac\frac \def\frac{\displaystyle\!frac}
\catcode`!=12

\let\oldsection\section
\renewcommand\section{\setcounter{equation}{0}\oldsection}

\allowdisplaybreaks
\def\pf{\it{Proof.}\rm\quad}

\def\R{\mathbb{R}}

\newtheorem{thm}{Theorem}[section]

\newtheorem{pro}{Proposition}[section]
\newtheorem{lem}{Lemma}[section]

\begin{document}
\title{\Large\bf
Asymptotic stability of strong  contact discontinuity  for full
compressible  Navier-Stokes equations  with initial boundary value problem    }
\author{{Tingting Zheng\footnote{Corresponding author email:
asting16@sohu.com(T.Zheng),}\footnote{
 This work was partially supported by
National Natural Science Foundation of China (Grant No. 11426062)
.}\footnote{This work was partially supported by the Fujian Province
Education Department
Project(JA12104).}\quad Yurui Lin}\\[2mm]
{\small  Computer and Message Science College, Fujian Agriculture and Forestry University,}\\[2mm]
{\small  Fuzhou 350001, P. R. China }}

\date{}

\maketitle

\noindent{\bf Abstract.} This paper is concerned with Dirichlet problem $u(0,t)=0$, $\theta(0,t)=\theta_-$
 for one-dimensional full compressible
Navier-Stokes equations  in
the half space $\R_+=(0,+\infty)$. Because the boundary decay rate is hard to control, stability of contact discontinuity result is very difficult. In this paper, we raise the decay rate and establish that for a certain class of large perturbation, the asymptotic stability result is
contact discontinuity. Also, we ask the strength of
contact discontinuity  not small.  The proofs are given by the elementary
energy method.
\\[2mm]
\section{Introduction} We consider
one-dimensional compressible viscous heat-conducting flow in the
half space $\R_+=[0,\infty)$, which is governed by the following
initial-boundary value problem in Eulerian coordinate
$(\tilde{x},t)$:
\begin{equation}
\left\{
\begin{array}{lll}
\tilde\rho_t+(\tilde\rho \tilde u)_{\tilde x}=0,\quad (\tilde x,t)\in\R_+\times\R_+,\\[2mm]
(\tilde\rho \tilde u)_t+(\tilde\rho \tilde u^2+\tilde p)_{\tilde{x}}=\mu \tilde u_{\tilde{x}\tilde{x}}, \\[2mm]
\left(\tilde\rho\left(\tilde e+\frac{\tilde
u^2}{2}\right)\right)_t+\left(\tilde\rho \tilde u\left(\tilde
e+\frac{\tilde u^2}{2}\right)+\tilde p\tilde
u\right)_{\tilde{x}}=\kappa\tilde
\theta_{\tilde{x}\tilde{x}}+(\mu\tilde u\tilde
u_{\tilde{x}})_{\tilde{x}},
\\[4mm]
(\tilde\rho,\tilde u,\tilde \theta)|_{\tilde{x}=0}=(\rho_-,0,\theta_-), \\[2mm]
(\tilde \rho,\tilde u,\tilde \theta)|_{t=0}=(\tilde\rho_0,\tilde
u_0,\tilde\theta_0)(\tilde x)\to(\rho_+,0,\theta_+)
\quad\mbox{as}\quad \tilde{x}\to\infty,
\end{array}\right.
\label{1.1}
\end{equation}
where $\tilde\rho$, $\tilde u$ and $\tilde\theta$ are the density,
the velocity and the absolute temperature, respectively, while
$\mu>0$ is the viscosity coefficient and $\kappa>0$ is the
heat-conductivity coefficients, respectively. It is assumed
throughout the paper that $\rho_\pm$ and $\theta_\pm$ are prescribed
positive constants . We shall focus our interest on the
 polytropic ideal gas  with
$|\theta_+-\theta_-|$ and  $|\rho_+-\rho_-|$ are general constants
    , so  the pressure $\tilde p=\tilde
p(\tilde\rho,\tilde\theta)$ and the internal energy $\tilde e=\tilde
e(\tilde\rho,\tilde\theta)$ are related by the second law of
thermodynamics:
\begin{equation}\label{1.2}
\tilde p=R\tilde\rho\tilde\theta,\quad \tilde
e=\frac{R}{\gamma-1}\tilde\theta+const.,
\end{equation}
where $\gamma>1$ is the adiabatic exponent and $R>0$ is the gas
constant.

The boundary condition $(\ref{1.1})_4$ implies that, through the
 boundary $\tilde
x=0$, the fluid with density $\rho_-$ flows into the region $\R_+$
at the speed $u=0$. So the initial-boundary value problem
(\ref{1.1}) is the so-called {\it impermeable wall} problem. In
terms of various boundary values,  Matsumura \cite{Ma2001}
classified all possible large-time behaviors of the solutions for
the one-dimensional (isentropic)compressible Navier-Stokes
equations.

To state our main results we
first transfer (\ref{1.1}) to the problem in the Lagrangian
coordinate and then make use of a coordinate transformation to
reduce the initial-boundary value problem (\ref{1.1}) into the
following form:
\begin{equation}
\left\{
\begin{array}{ll}
v_t-u_x=0,\quad (x,t)\in\R_+\times\R_+,\\[2mm]
u_t+\left(\frac{R\theta}{v}\right)_x=\mu\left(\frac{u_x}{v}\right)_x,
\\[4mm]
\frac{R}{\gamma-1}\theta_t+R\frac{\theta}{v}u_x=\kappa\left(\frac{\theta_x}{v}\right)_x+\mu\frac{u_x^2}{v},\\[4mm]
(v,u,\theta)|_{x=0}=(v_-,0,\theta_-),\quad t>0,\\[2mm]
(v,u,\theta)|_{t=0}=(v_0,u_0,\theta_0)\to(v_+,0,\theta_+)\quad\mbox{as}\quad
x\to\infty,
\end{array}
\right.\label{1.3}
\end{equation}
where $v_\pm$ and $\theta_\pm$ are given positive constants, and
 $v_0,\ \theta_0>0$. In fact $v=1/\rho(x,t),\
u=u(x,t),\ \theta=\theta(x,t)$ and $R\theta/v=p(v,\theta)$ are the
specific volume, velocity , temperature and pressure as in
(\ref{1.1}).

There have been a lot of works on the asymptotic behaviors of
solutions to the initial-boundary value (or Cauchy) problem for the
Navier-Stokes equations toward  basic waves or their viscous
versions, see, for example, [3--26] and the reference therein. In terms of various boundary values, Matsumura and Nishihara \cite{MN2001} classified all possible large-time behaviors of the solutions for one-dimensional (isentropic) compressible Navier-Stokes equations. In the case where $u(0,t)=0$, the problem is called $impermeable\  wall$ problem.
Inflow problem is one of Dirichlet problems with $u(0,t)<0$. Matsumura and Nishihara \cite{MN2001} have obtained the stability theorems on both the boundary layer solution and the superposition of a boundary-layer solution and a rarefaction wave for inflow problem.  Due to Huang et al. \cite{HMS2003} in which the
asymptotic stability on both the viscous shock wave and the
superposition of a viscous shock wave and a boundary-layer solution
are studied. However, the problem of stability of contact
discontinuities are associated with linear degenerate fields and previous results are
less stable than the nonlinear waves for the inviscid system (Euler
equations).  It was observed in \cite{X,LX}, where the metastability
of contact waves was studied for viscous conservation laws with
artificial viscosity, that the contact discontinuity cannot be the
asymptotic state for the viscous system, and a diffusive wave, which
approximated the contact discontinuity on any finite time interval,
actually dominates the large-time behavior of solutions. The
nonlinear stability of weak contact discontinuity for the (full)
compressible Navier-Stokes equations was then investigated in
\cite{HMS,HZ}
 for the free  boundary value problem , \cite{HMX,HXY} for the Cauchy problem and \cite{TJJ2011,WQ} for the inflow problem.
In \cite{HMX}£¬ they point out  because of the decay rate problem, viscous contact discontinuity for Dirichlet
impermeable wall problem becomes difficulty.

Except Dirichlet
impermeable wall problem, recently, some problems are call stability of strong viscous
waves(see \cite{NYZ}--\cite{Z2013}, \cite{T2015}). These stability
results (to Cauchy problem or free boundary problem) are shown with some special conditions. Especially,  zero
dissipation result is shown in\cite{MaSX} and $\gamma\to 1$ in
\cite{NYZ} or \cite{HH}. Base on small oscillation, initial
smallness perturbation or zero dissipation (and so on),
Navier-Stokes equations stability results can be obtained with some
skills.

The main purpose of this paper is to improve the previous studies for Dirichlet problem . Base on the new estimates on
the heat kernel in \cite{TZ2012}, the decay rate can be promote. Then, we justify that the solution $(v,u,\theta)$ of the Navier-Stokes system (\ref{1.3}) tends asymptotically to contact discontinuity which is the Riemann solution  (\ref{1.5}) in $\R_+$. For the strength of the contact discontinuity is not small, we call it strong contact discontinuity .
Also,  we  get rid off  smallness $\|(\varphi_{0x},\psi_{0x},\zeta_{0x})\|_{L^2}$
and $\gamma\to 1$, i.e.,  for the Direchlet  problem (\ref{1.3}), it is possible to be resolved and
stable and the solution approximate the strong contact
discontinuity .

$ \mathbf{Notation.}$ Throughout this paper, we shall denote
$H^l(\R_+)$ the usual $l-th$ order Sobolev space with the norm
$$\|f\|_{l}=\big(\sum_{j=0}^l\|\partial_x^j f\|^2\big)^{1/2},\ \ \|\cdot\|:=\|\cdot\|_{L^2(\R_+)}.$$
For simplicity, we also use $C$ or $C_i$ ($i=1,2,3.....$) to denote
the various positive generic constants.  $C(\delta_0)$ stands for
suitably small constant about $\delta_0$ and $C_v=\frac{R}{\gamma-1}$ . $\epsilon$ and
$\epsilon_i$($i=1,2,3.....$) stand for suitably small positive
constant in Cauchy-Schwarz inequality and
$\partial_x^i=\frac{\partial^i}{\partial x^i}$.

\section{Reformulation and main result}
As shown in previous studies, the
asymptotic behavior is well characterized by the solutions to the
corresponding Riemann problem for the hyperbolic part of (\ref{1.3})
(that is, Euler system):

\begin{equation}
\left\{
\begin{array}{ll}
v_t-u_x=0,\\[2mm]
u_t+p(v,\theta)_x=0,\\[2mm]
\frac{R}{\gamma-1}\theta_t+R\frac{\theta}{v}u_x=0,\\[4mm]
(v,u,\theta)(x,0)=(v_-,0,\theta_-)\quad\mbox{if}\quad x<0,\\[2mm]
(v,u,\theta)(x,0)=(v_+,0,\theta_+)\quad\mbox{if}\quad x>0,
\end{array}
\right.\label{1.4}
\end{equation}
The Riemann problem of system (\ref{1.4}) admits a contact
discontinuity

\begin{equation}
\left(\overline{V},\overline{U},\overline{\Theta}\right) =\left\{
\begin{array}{ll}
(v_-,0,\theta_-),\quad x<0,\\[2mm]
(v_+,0,\theta_+),\quad x>0,
\end{array}
\right. \label{1.5}
\end{equation}
provided that
\begin{equation}
p_-=R\frac{\theta_-}{v_-}=p_+=R\frac{\theta_+}{v_+}.\label{1.6}
\end{equation}
As that in \cite{HMS} we conjecture a pair of  $(V,U,\Theta)(x,t)$  is as follows

\begin{equation}
P(V,\Theta)=R\frac{\Theta}{V}=p_+,\;\;\;U(x,t)=\frac{\kappa(\gamma-1)\Theta_x}{\gamma
R\Theta},\label{1.7}
\end{equation}
and
\begin{equation}
\left\{ \begin{array}{ll}
\Theta_t=a(\ln\Theta)_{xx},\quad a=\frac{\kappa p_+(\gamma-1)}{\gamma R^2}>0,\\[2mm]
\Theta(0,t)=\theta_-,\\[2mm]
 \Theta(x,0)=\Theta_{0}(x)\to \theta_+,\ \ as\ \ x\to +\infty,
\end{array}\right.\label{1.8}
\end{equation}
with
$\Theta_0(x)=\left(\frac{2}{\sqrt{\pi}}(\theta_+^{1/\delta_0}-\theta_-^{1/\delta_0})\int_0^{\ln(x+\sqrt{1+x^2})}\exp\{-y^2\}dy
 +\theta_-^{1/\delta_0}\right)^{\delta_0}.$ Here $\delta_0$ is a
 suitably small constant and $1/\delta_0$ is an integer.

 According to the smallness $\delta_0$, we can find $\Theta_0(x)$ is
 nearly like a function $f(x)=\left\{
                                        \begin{array}{ll}
                                          \theta_+, & \hbox{$x>0$;} \\
                                          \theta_-, & \hbox{$x=0$.}
                                        \end{array}
                                      \right.
 $
That means $\Theta_0(x)$ satisfying the following properties.
\begin{lem}\label{lem2.1}
\begin{eqnarray}\label{ca1.12}
&& \|\Theta_{0x}\|_{L^1}\leq C, \ |\Theta_{0x}|\leq C\delta_0,\ \ |\Theta_{0xx}|\leq C\delta_0,\ \ \|\Theta_{0x}\|^2\leq C\delta_0^2,\nonumber\\
&&\|\Theta_{0xx}\|^2\leq
C\delta^2_0,\ \ \|\Theta_{0xxx}\|^2\leq C,\ \ \|\Theta_0-\theta_+\|_{L^1}\leq C.
\end{eqnarray}
\end{lem}

\pf In fact, if  $K(x)=\ln(x+\sqrt{1+x^2})$, we can get
\begin{eqnarray}\label{2.2}
\int_{\R_+}\exp\{-K^2(x)\}dx&=&\int_{\R_+}\exp\{-K^2(x)\}\frac{\sqrt{1+x^2}}{\sqrt{1+x^2}}dx\nonumber\\
&&=\int_{\R_+}\exp\{-K^2(x)\}\sqrt{1+x^2}d K(x)\nonumber\\
&&\leq\int_{\R_+}\exp\{-K^2(x)\}\exp\{K(x)\}d K(x)\nonumber\\
&&\leq\int_{\R_+}\exp\{-K^2(x)+K(x)\}d K(x)\leq C.
\end{eqnarray}Set
\begin{equation}\label{2.3}H(x)=\Theta_0^{1/\delta_0}=\frac{\theta_+^{1/\delta_0}+\theta_-^{1/\delta_0}}{2}
+\frac{\theta_+^{1/\delta_0}-\theta_-^{1/\delta_0}}{\sqrt{\pi}}\int_0^{K(x)}\exp\{-x^2\}dx,\end{equation}
from $K_x=(1+x^2)^{-1/2}$ and
\begin{equation}\label{2.4}\Theta_{0x}=\delta_0H^{\delta_0-1}(x)H_x(x)=\delta_0H(x)^{\delta_0-1}\frac{\theta_+^{1/\delta_0}-\theta_-^{1/\delta_0}}{\sqrt{\pi}}K_x(x)
\exp\{-K^2(x)\}>0,\ \ x\in\R_+,\end{equation} we can get
$\|\Theta_{0x}\|_{L^1(\R_+)}<C$.

 When $x>0$, $K(x)>0$, we can know
\begin{equation}\label{2.5}\frac{\theta_+^{1/\delta_0}-\theta_-^{1/\delta_0}}{H(x)}\leq
C\frac{\theta_+^{1/\delta_0}-\theta_-^{1/\delta_0}}{\theta_+^{1/\delta_0}+\theta_-^{1/\delta_0}}\leq
C .\end{equation} Because
\begin{equation}\label{2.8a}\theta_-^{1/\delta_0}\leq H(x)\leq
\theta_+^{1/\delta_0},\end{equation} from (\ref{2.4}) and
$K_x(x)=\frac{1}{\sqrt{1+x^2}}$,
 we can get $
|\Theta_{0x}|\leq C\delta_0$.  Also from (\ref{2.2}),(\ref{2.3}),
(\ref{2.4}) and (\ref{2.5}) we can get $\|\Theta_{0x}\|\leq
C\delta_0$. Similar as above estimates, it is easy to  check
$\|\Theta_{0xx}\|\leq C\delta_0^2$ and $\|\Theta_{0xxx}\|\leq C$.

When $\delta_0=\frac{1}{2k+1}$, $k\in \{1,2,3...\}$ is a suitably
large constant, from the  equality
$a^n-b^n=(a-b)\sum_{i=0}^{n-1}a^{n-1-i}b^i$, $\forall a>0,\ b>0,\
n\in\{1,2,3...\}$ and (\ref{2.2}),(\ref{2.4}), (\ref{2.5}), we can
get that
\begin{eqnarray}\label{2.6}&&\int_{\R_+}|\Theta_0-\theta_+|dx=\int_{\R_+}|H^{\delta_0}-\theta_+|dx\nonumber\\
&&=\int_{\R_+}\frac{|H-\theta_+^{1/\delta_0}|}{\sum_{i=0}^{2k}H^{(2k-i)/(2k+1)}\theta_+^{i/(2k+1)}}dx\nonumber\\
&&\leq
C\int_{\R_+}\frac{(\theta_+^{1/\delta_0}-\theta_-^{1/\delta_0})\exp\{-CK^2(x)\}}{\sum_{i=0}^{2k}H^{(2k-i)/(2k+1)}\theta_+^{i/(2k+1)}}dx\nonumber\\
&&\leq
C(\theta_+^{1/\delta_0}-\theta_-^{1/\delta_0})\sup_{x\in\R_+}H^{1/(2k+1)-1}\leq
C(\theta_++\theta_-).\end{eqnarray} $\Box$

 In summary we have constructed a pair of functions
$(V,U,\Theta)$ such that
\begin{equation}
\left\{
\begin{array}{ll}
R\frac{\Theta}{V}=p_+,\\[4mm]
V_t=U_x,\\[2mm]
U_t+P(V,\Theta)_x=\mu\left(\frac{U_x}{V}\right)_x+F,\\[4mm]
\frac{R}{\gamma-1}\Theta_t+R\frac{\Theta}{V}U_x=\kappa\left(\frac{\Theta_x}{V}\right)_x+\mu\frac{U_x^2}{V}+G,\\[4mm]
(V,U,\Theta)(0,t)=(v_-,\frac{\kappa(\gamma-1)}{\gamma
R}\frac{\Theta_{x}}{\Theta}|_{x=0},\theta_-),\\[4mm]
(V,U,\Theta)(x,0)=(V_0,U_0,\Theta_0)=(\frac{R}{p_+}\Theta_0,\frac{\kappa(\gamma-1)}{\gamma
R}\frac{\Theta_{0x}}{\Theta_0},\Theta_0)\to (v_+,0,\theta_+),\ as\ \
x\to+\infty,
\end{array}\right.\label{1.9}
\end{equation}
where
\begin{eqnarray}G(x,t)&=&-\mu\frac{U_x^2}{V}=O((\ln\Theta)^2_{xx}),\nonumber\\
F(x,t)&=&\frac{\kappa(\gamma-1)}{\gamma
R}\left\{(\ln\Theta)_{xt}-\mu\left(\frac{(\ln\Theta)_{xx}}{V}\right)_x\right\}\nonumber\\
&=&\frac{\kappa a(\gamma-1)-\mu p_+\gamma}{R
\gamma}\left(\frac{(\ln\Theta)_{xx}}{\Theta}\right)_x.\label{1.10}\end{eqnarray}

Denote
\begin{eqnarray}
&&\varphi(x,t)=v(x,t)-V(x,t),\nonumber\\
&&\psi(x,t)=u(x,t)-U(x,t),\nonumber\\
&&\zeta(x,t)=\theta(x,t)-\Theta(x,t).\label{1.11}
\end{eqnarray}
Combining (\ref{1.9}) and (\ref{1.3}), the original problem can be
reformulated as
\begin{equation}\left\{
    \begin{array}{lll}
      \varphi_t=\psi_x, &  \\
      \psi_t-(\frac{R\Theta}{vV}\varphi)_x+(\frac{R\zeta}{v})_x=-\mu(\frac{U_x}{vV}\varphi)_x+\mu(\frac{\psi_x}{v})_x-F, & \\
      \frac{R}{\gamma-1}\zeta_t+\frac{R\theta}{v}(\psi_x+U_x)-\frac{R\Theta}{V}U_x
=\kappa(\frac{\zeta_x}{v})_x-\kappa(\frac{\Theta_x\varphi}{vV})_x+\mu(\frac{{u_x}^2}{v}-\frac{{U_x}^2}{V})-G,
&\\
(\varphi,\psi,\zeta)(0,t)=(0,u_b-U(0,t),0)=(0,-\kappa(\gamma-1)\Theta_x(0,t)/(\gamma R\theta_-),0),&\\
(\varphi,\psi,\zeta)(x,0)=(\varphi_0,\psi_0,\zeta_0)=(v_0-V_0,u_0-U_0,\theta_0-\Theta_0),&
\end{array}
  \right.
\label{1.12}\end{equation}

Under the above preparation in hand,  we assume throughout of this
section that
$$(\varphi_0,\zeta_0)(x)\in H_0^1(0,\infty),\ \ \psi_0(x)\in H^1(0,\infty).$$
Moreover, for an interval $I\in [0,\infty)$ , we define the function
space
$$X(I)=\left\{(\varphi,\psi,\zeta)\in C(I,H^1)|\varphi_x\in
L^2(I;L^2), (\psi_x,\zeta_x)\in L^2(I;H^1)\right\}.$$
 Our main results of this paper now reads as
follows.

\begin{thm}\label{thm1.1}If
  $(v_0-v_+,u_0,\theta_0-\theta_+)\in H^2(\R_+)\cap L^1(\R_+)$, $\|(v_0-\overline{V},u_0-\overline{U},\theta_0-\overline{\Theta})\|$ suitably small,
$\frac{v_-}{\theta_-}=\frac{v_+}{\theta_+}$  and
$|\theta_+-\theta_-|$ not small,
  (\ref{1.3}) has a global solution $(v,u,\theta)$  satisfying
$(\varphi,\psi,\zeta)\in X([0,\infty))$, and when $t\to\infty$,
$$\|(v-\overline{V},u-\overline{U},\theta-\overline{U})\|_{L^{\infty}(\R_+)}\to (0,0,0).$$
\end{thm}

\section{Preliminary}

In this section, to study the asymptotic behavior of the solution to
the Cauchy problem (\ref{1.3}), we provide some preliminary lemmas
and list the a priori estimate that are important for the proof of
Theorem \ref{thm1.1}.

\begin{lem}\label{lem2.2}If $\delta_0$ and $\Theta_0$ satisfying the
condition in Theorem \ref{thm1.1} and
\begin{eqnarray*}&&\theta_2(x,t)=\int_{0}^{+\infty}(4\pi
at)^{-1/2}(\Theta_0(h)-\theta_-)\left\{\exp\{-\frac{(h-x)^2}{4at}\}-\exp\{-\frac{(h+x)^2}{4at}\}\right\}\
dh+\theta_-,\end{eqnarray*}
 we can get
 \begin{eqnarray}
 &&\theta_{2t}=a\theta_{2xx};\nonumber\\
 &&\theta_2(0,t)=\theta_-;\nonumber\\
 &&\theta_2(x,0)=\theta_{20}(x)=\left\{
                      \begin{array}{ll}
                        \Theta_{0}(x)\to \theta_+ , & \hbox{$x>0$;} \\
                       -\Theta_{0}(-x)+2\theta_-\to 2\theta_--\theta_+ , & \hbox{$ x\leq 0$,}
                      \end{array}
                    \right.\label{2.2}
 \end{eqnarray}
 and\begin{eqnarray}
 &&\int_0^t\|\theta_{2x}\|^2dt\leq C(1+t)^{1/3},\label{2.3}\end{eqnarray}

\end{lem}
\pf Because $\theta_2(x,t)$ can be rewrite to
$$\theta_2(x,t)=\int_{-\infty}^{+\infty}(4\pi at)^{-1/2}\theta_{20}(h)\exp\{-\frac{(x-h)^2}{4at}\}dh,$$
and $\theta_{20}(x)\in  C^1(\R)$, we find that $\theta_2(x,t)$ is a
fundamental solution of (\ref{2.2}), it is easy to check $\lim_{t\to
0}\theta_2(x,t)=\theta_{20}(x)$, so we finish (\ref{2.2}).

Because \begin{eqnarray}\label{2.5} \theta_{2x}
&=&\int_0^{+\infty}(4\pi
at)^{-1/2}\left(\Theta_0(z)-\theta_-\right)\Big\{\exp\{-\frac{(z-x)^2}{4at}\}\frac{z-x}{2at}+\exp\{-\frac{(z+x)^2}{4at}\}\frac{z+x}{2at}\Big\}dz\nonumber\\
 &=&\int_0^{\infty}(4\pi
at)^{-1/2}\Theta_{0z}(z)\Big\{\exp\{\frac{-(z-x)^2}{4at}\}-\exp\{\frac{-(z+x)^2}{4at}\}\Big\}dz\nonumber\\
&=&\int_0^{+\infty}(4\pi
at)^{-1/2}\left(\Theta_0(z)-\theta_++\theta_+-\Theta_0(x)\right)\nonumber\\
&&\quad\times\Big\{\exp\{-\frac{(z-x)^2}{4at}\}\frac{z-x}{2at}+\exp\{-\frac{(z+x)^2}{4at}\}\frac{z+x}{2at}\Big\}dz,\end{eqnarray}

 By using H$\ddot{o}$lder inequality , Fubini Theorem and
$\|\Theta_{0}-\theta_+\|_{L^1(\R_+)}<C$ , we can get from
(\ref{2.5}) that
\begin{eqnarray*}\int_0^t\int_0^{\infty}\theta_{2x}^2dxdt&\leq&
C\int_0^t\int_0^{\infty}(4\pi
at)^{-1}\left\{\int_0^{\infty}\Theta_{0z}\left(\exp\{-\frac{(z-x)^2}{4at}\}-\exp\{-\frac{(z+x)^2}{4at}\}\right) dz\right\}^2dxdt\\
&\leq&C\int_0^t\int_0^{\infty}(4\pi
at)^{-1}\int_0^{\infty}|\Theta_{0z}|\Big\{\exp\{-\frac{(z-x)^2}{4at}\}
+\exp\{-\frac{(z+x)^2}{4at}\} \Big\}dz\\
&&\quad\times\int_0^{+\infty}(4\pi
at)^{-1/2}\left|\Theta_0(z)-\theta_++\theta_+-\Theta_0(x)\right|\Big\{\exp\{-\frac{(z-x)^2}{4at}\}\frac{|z-x|}{2at}\nonumber\\
&&\quad\quad+\exp\{-\frac{(z+x)^2}{4at}\}\frac{|z+x|}{2at}\Big\}dzdxdt\\
&&\quad +C\int_0^1\int_0^{\infty}(4\pi
at)^{-1/2}\int_0^{\infty}|\Theta_{0z}|\Big\{\exp\{-\frac{(z-x)^2}{4at}\}
+\exp\{-\frac{(z+x)^2}{4at}\}\Big\}
dzdxdt\\
&\leq&C\int_1^t\int_0^{\infty}(4\pi at)^{-1}\int_0^{+\infty}(4\pi
at)^{-1/2}\left|\Theta_0(z)-\theta_+\right|\\
&&\times\Big\{\exp\{-\frac{(z-x)^2}{4at}\}\frac{|z-x|}{2at}+\exp\{-\frac{(z+x)^2}{4at}\}\frac{|z+x|}{2at}\Big\}dxdzdt\nonumber\\
&&+C\int_1^t\int_0^{\infty}(4\pi at)^{-1}\int_0^{+\infty}(4\pi
at)^{-1/2}\left|\theta_+-\Theta_0(x)\right|\\
&&\times\Big\{\exp\{-\frac{(z-x)^2}{4at}\}\frac{|z-x|}{2at}+\exp\{-\frac{(z+x)^2}{4at}\}\frac{|z+x|}{2at}\Big\}dzdxdt+C\\
 &\leq &C\|\Theta_0-\theta_+\|_{L^1(\R_+)}\ln(1+t)+C\leq C(1+t)^{1/3}.\end{eqnarray*}

So we finish this lemma.$\Box$

From (\ref{1.7}) and (\ref{1.9}), we obtain
\begin{eqnarray}\label{2.15a}
(|V_x|+|U|)\leq C|\Theta_x|,\   |\Theta_x|^2\leq
C\|(\ln\Theta)_x\|\|(\ln\Theta)_{xx}\|,\ |U_x|^2\leq
C\|(\ln\Theta)_{xx}\|\|(\ln\Theta)_{xxx}\|.
\end{eqnarray} According to the definition of $(\overline{V},\overline{U},\overline{\Theta})$ in (\ref{1.5}),
when the time $t\to\infty$ and $x\in\R_+$, it is easily check that
$(V,U,\Theta)$ is nearly close to the contact discontinuity $(\overline{V},\overline{U},\overline{\Theta})$ by the following lemma.

\begin{lem}\label{lem2.3}
There exist a positive constant  $C$ such that
\begin{eqnarray}&&\|(\ln\Theta)_x\|^2+a\int_0^t\ \|(\ln\Theta)_{xx}\|^2\ dt \leq
C\delta_0^2.\label{2.7}\\
&&_{(see (\ref{2.17})-(\ref{2.19}))}\nonumber\\
 &&\|\Theta-\theta_2\|^2+\int_0^t\ \|(\ln\Theta)_{x}\|^2\ dt\leq
C(1+t)^{1/3}.\label{2.8}\\
&&_{{(see (\ref{2.14})-(\ref{2.15}))}}\nonumber\\
 && \|(\ln\Theta)_x\|^2\leq
C(1+t)^{-2/3}.\label{2.9}\\&&_{{(see
(\ref{2.20})-(\ref{2.23}))}}\nonumber\\
&&\|(\ln\Theta)_{xx}\|^2\leq C(1+t)^{-5/3}.\label{2.10}\\&&_{(see
(\ref{2.24})-(\ref{2.28}))}\nonumber\\
&&\|(\ln\Theta)_{xx}\|^2(1+t)+\int_0^t\|\partial^3_x\ln\Theta\|^2(1+t)\
dt\leq
C\delta_0^2.\label{2.11}\\&&_{(see(\ref{2.29}))}\nonumber\\
&& \|\partial^3_x\ln\Theta\|^2\leq
C(1+t)^{-8/3}.\label{2.12}\\&&_{(see
(\ref{2.30})-(\ref{2.32}))}\nonumber\\
&&\|\Theta-\theta_+\|^2_{L^{\infty}}\leq
C\delta_0^{1/4}(1+t)^{-1/24}.\label{2.13b}\\ &&_{(see
(\ref{3.29a})-(\ref{3.31a}))}\nonumber
\end{eqnarray}
\end{lem}
\pf From (\ref{1.8}) we know
\begin{equation}(\ln\Theta)_t=a\frac{(\ln\Theta)_{xx}}{\Theta},\label{2.16}\end{equation}
both side of it multiply by $(\ln\Theta)_{xx}$ and integrate in
$\R_+\times (0,t)$ we can get

\begin{eqnarray}&&\|(\ln\Theta)_x\|^2+a\int_0^t\ \|(\ln\Theta)_{xx}\|^2\
dt\nonumber\\
&&\leq
\|(\ln\Theta_{0})_x\|^2+\int_0^t(\ln\Theta)_t(\ln\Theta)_x\big|_0^{\infty}\
dt.\label{2.17}\end{eqnarray}

Then from (\ref{ca1.12}) and (\ref{2.17}) we can get
\begin{equation}\|(\ln\Theta)_x\|^2+a\int_0^t\ \|(\ln\Theta)_{xx}\|^2\ dt\leq
C\delta_0.\label{2.19}\end{equation}

Then,  if $\int_0^t\int_{\R_+}
\left((\ref{1.8})_1-(\ref{2.2})_1\right)\times(\Theta-\theta_2)dxdt$
  combine with Cauchy-Schwarz inequality we can get
\begin{eqnarray}&&\|\Theta-\theta_2\|^2+\int_0^t\|(\ln\Theta)_x\|^2\
dt\leq C\int_0^t\|{\theta_2}_x\|^2dt.\label{2.14}\end{eqnarray} Use
(\ref{2.3}) to (\ref{2.14}) we can get
\begin{equation}\|\Theta-\theta_2\|^2+\int_0^t\ \|(\ln\Theta)_{x}\|^2\ dt\leq
C(1+t)^{1/3}. \label{2.15}\end{equation} That is (\ref{2.8}).

 Next, from $$\int_0^t\int_{\R_+}(\ref{1.8})_1 \times
\Theta^{-1}(\ln\Theta)_{xx}(1+t)dxdt,$$  we can get
\begin{eqnarray}&&\int_0^t(1+t)\left((\ln\Theta)_t(\ln\Theta)_x\right)(0,t)\ dt\nonumber\\
&&=a\int_0^t\int_0^{\infty}\ \frac{(\ln\Theta)^2_{xx}}{\Theta}(1+t)\
dxdt+\int_0^t\int_0^{\infty}\ \left((\ln\Theta)^2_x\right)_t(1+t)\
dxdt.\label{2.20}\end{eqnarray} Because
\begin{equation}\label{2.21}\int_0^t(1+t)(\ln\Theta)_t(\ln\Theta)_x(0,t)\
dt=0,\end{equation} we can get
\begin{eqnarray}
&&(1+t)\|(\ln\Theta)_x\|^2+\int_0^t\
\int_0^{\infty}(1+t)(\ln\Theta)^2_{xx}\ dx\ dt\nonumber\\
 &&\leq C\|\Theta_{0x}\|^2+\int_0^t\ \int_0^{\infty}\ (\ln\Theta)^2_x\ dx\ dt.\label{2.22}\end{eqnarray}
Combine with (\ref{2.15}) we can get
\begin{eqnarray}
&&(1+t)\|(\ln\Theta)_x\|^2+\int_0^t\
\int_0^{\infty}(1+t)(\ln\Theta)^2_{xx}\ dx\ dt\nonumber\\
 &&\leq C(1+t)^{1/3}.\label{2.23}\end{eqnarray}
That means $\|(\ln\Theta)_x\|^2\leq C(1+t)^{-2/3}$, which
is(\ref{2.9}).

Again from (\ref{1.8})$_1$ we can get
\begin{equation}(\ln\Theta)_{xt}=a\left(\frac{(\ln\Theta)_{xx}}{\Theta}\right)_x.\label{2.24}
\end{equation}
Both side of (\ref{2.24})multiply $\partial^3_x\ln\Theta$ and get
\begin{equation}\left((\ln\Theta)_{xt}\partial^2_x(\ln\Theta)\right)_x-1/2(\partial^2_x\ln\Theta)_t
=a\left(\frac{(\ln\Theta)_{xx}}{\Theta}\right)_x\partial^3_x(\ln\Theta).\label{2.25}\end{equation}
Because
\begin{eqnarray*}&&\left((\ln\Theta)_{xt}\partial^2_x(\ln\Theta)\right)_x(1+t)^2\\
&&=\left((\ln\Theta)_{xt}(\ln\Theta)_{xx}\right)_x(1+t)^2\\
&&=a^{-1}\left((\ln\Theta)_{xt}\Theta_t\right)_x(1+t)^2,\end{eqnarray*}
both side of (\ref{2.25}) multiply $(1+t)^2$ , then integrate in
$\R_+\times(0,t)$ and combine with $\Theta_t(0,t)=0$,
$\Theta_t(\infty,t)=0$, $\Theta_x(\infty,t)=0$ and Cauchy-Schwarz
inequality to get for some small $\epsilon>0$ we have
\begin{eqnarray}
&&0\geq a\int_0^t\ \int_0^{\infty}\
\frac{(\ln\Theta)_{xxx}^2}{\Theta}(1+t)^2\ dx\
dt\nonumber\\
&&\quad-\epsilon\int_0^t\ \int_0^{\infty}\
(1+t)^2(\ln\Theta)_{xxx}^2\ dx\ dt-C\epsilon^{-1}a\int_0^t\
\int_0^{\infty}\ (1+t)^2(\ln\Theta)_{xx}^2(\ln\Theta)_x^2\ dx\
dt\nonumber\\
&&\quad+1/2\|(\ln\Theta)_{xx}\|^2(1+t)^2-1/2\|(\ln\Theta_0)_{xx}\|^2-\int_0^t\
\|(\ln\Theta)_{xx}\|^2(1+t)\ dx\nonumber\\
&&\geq Ca\int_0^t\ \int_0^{\infty}\
\frac{(\ln\Theta)_{xxx}^2}{\Theta}(1+t)^2\ dx\
dt\nonumber\\
&&\quad-C\epsilon^{-1}a\int_0^t\ \int_0^{\infty}\
(1+t)^2\|(\ln\Theta)_{xx}\|\|(\ln\Theta)_{xxx}\|(\ln\Theta)_x^2\ dx\
dt\nonumber\\
&&\quad+1/2\|(\ln\Theta)_{xx}\|^2(1+t)^2-1/2\|(\ln\Theta_0)_{xx}\|^2-\int_0^t\
\|(\ln\Theta)_{xx}\|^2(1+t)\ dx.\label{2.26}
\end{eqnarray}
 Take (\ref{2.23})
 into (\ref{2.26}) we can get

\begin{eqnarray}&&\|(\ln\Theta)_{xx}\|^2(1+t)^2+\int_0^t\ \int_0^{\infty}\
(1+t)^2(\ln\Theta)_{xxx}^2\ dx\ dt\nonumber\\
&& \leq C(1+t)^{1/3},\label{2.27}\end{eqnarray} which also means

\begin{equation}\|(\ln\Theta)_{xx}\|^2\leq
C(1+t)^{-5/3},\label{2.28}\end{equation} so we finish (\ref{2.10}).

If both side of (\ref{2.25}) multiply by $(1+t)$, similar as the
proof of (\ref{2.27}),  when combine with (\ref{2.19}) we can get
\begin{equation}\|(\ln\Theta)_{xx}\|^2(1+t)+\int_0^t\ \int_0^{\infty}\ (1+t)(\partial^3_x\ln\Theta)^2\ dx\
dt\leq C\delta_0^2,\label{2.29}\end{equation} which means
(\ref{2.11}).

From (\ref{2.24}) we can get
\begin{equation}
\partial_t(\ln\Theta)_{xx}=a\partial^2_x\left(\frac{(\ln\Theta)_{xx}}{\Theta}\right).\label{2.30}
\end{equation}
Because
\begin{eqnarray*}
\left((\ln\Theta)_{xxt}(\ln\Theta)_{xxx}\right)_x=\left(a^{-1}\Theta_{tt}(\ln\Theta)_{xxx}\right)_x
,\end{eqnarray*}
when both side of (\ref{2.30}) multiply
$(\partial^4_x\ln\Theta)(1+\tau)^3$, then integrate in
$\R_+\times(0,t)$, we can get

\begin{eqnarray}
&&\int_0^t\int_0^{\infty}\left(a^{-1}\Theta_{tt}(\ln\Theta)_{xxx}\right)_x(1+\tau)^3dxd\tau\nonumber\\
&&=\int_0^t\int_0^{\infty}a\partial^2_x\left(\frac{(\ln\Theta)_{xx}}{\Theta}\right)\partial^4_x\ln\Theta(1+\tau)^3dxd\tau\nonumber\\
&&\quad+\int_0^t\int_0^{\infty}\frac{1}{2}\left((\partial^3_x\ln\Theta)^2\right)_t(1+\tau)^3dxd\tau.\label{2.31}
\end{eqnarray}

So from (\ref{2.23}) and (\ref{2.27}), we can get that for a small
$\epsilon>0$, (\ref{2.31}) can be change to
\begin{eqnarray*}
&&\|\partial^3_x\ln\Theta\|^2(1+t)^3+C\int_0^t\
(1+\tau)^3\|\partial^4_x\ln\Theta\|^2\ d\tau\\
&&\leq C+C\int_0^t\ \int_0^{\infty}\
(\partial^3_x\ln\Theta)^2(\ln\Theta)_x^2(1+\tau)^3\ dx\
d\tau+C\int_0^t\
\int_0^{\infty}\ (\ln\Theta)_{xx}^4(1+\tau)^3\ dx\ d\tau\\
&&\quad+C\int_0^t\ \int_0^{\infty}\
(\partial^2_x\ln\Theta)^2(\ln\Theta)_x^4(1+\tau)^3\ dx\
d\tau+C\int_0^t\
\int_0^{\infty}\ (\partial_x^3\ln\Theta)^2(1+\tau)^2\ dx\ d\tau\\
&&\leq C\int_0^t\
\|(\ln\Theta)_x\|^2\|\partial^3_x\ln\Theta\|\|\partial^4_x\ln\Theta\|(1+\tau)^3\
d\tau+C\int_0^t\
\|(\ln\Theta)_{xx}\|^3\|\partial^3_x\ln\Theta\|(1+\tau)^3\ d\tau\\
&&\quad+\int_0^t\
\|(\ln\Theta)_{xx}\|^4\|(\ln\Theta)_x\|^2(1+\tau)^3\
d\tau+C(1+t)^{1/3}\\
&&\leq \epsilon\int_0^t\ \|\partial^4_x\ln\Theta\|^2(1+\tau)^3\
d\tau+C\epsilon^{-1}\int_0^t\ \|\partial^3_x\ln\Theta\|^2(1+\tau)^2\ d\tau\\
&&\quad+C\epsilon^{-1} \int_0^t\
\|\partial^2_x\ln\Theta\|^2(1+\tau)\ d\tau+C(1+t)^{1/3}.
\end{eqnarray*}
Again using (\ref{2.23}) and (\ref{2.27}) we can get
\begin{equation}\|\partial^3_x\ln\Theta\|^2(1+t)^3+\int_0^t\
(1+\tau)^3\|\partial^4_x\ln\Theta\|^2\ d\tau\leq
C(1+t)^{1/3}.\label{2.32}\end{equation} This means (\ref{2.12})
finished.

Similar as above, when both side of (\ref{2.30}) multiply
$(\partial^4_x\ln\Theta)(1+\tau)^2$, then integrate in
$\R_+\times(0,t)$ we can get
\begin{equation}\|\partial^3_x\ln\Theta\|^2(1+t)^2+\int_0^t\
(1+\tau)^2\|\partial^4_x\ln\Theta\|^2\ d\tau\leq
C.\label{2.32b}\end{equation}

From (\ref{1.7}) we can
get\begin{equation*}(\Theta-\Theta_0)_t(\Theta-\Theta_0)=a\left((\ln\Theta)_x(\Theta-\Theta_0)\right)_x-a(\ln\Theta)_x(\Theta-\Theta_0)_x.\end{equation*}
When integrate both sides of it integrate in $\R_+\times[0,t]$, we can get
\begin{eqnarray}\|\Theta-\Theta_0\|^2\leq
C\|\Theta_{0x}\|_{L^1}\int_0^t\|\Theta_x\|_{L^{\infty}}d\tau\leq
C\|\Theta_{0x}\|_{L^1}\int_0^t\|\Theta_x\|^{1/2}\|\Theta_{xx}\|^{1/2}d\tau.\label{3.29a}\end{eqnarray}
From Lemma \ref{lem2.1},(\ref{3.29a}), (\ref{2.9}) and (\ref{2.10})
we can obtain
\begin{equation}\|\Theta-\theta_+\|^2\leq C(1+t)^{5/12}+\|\theta_+-\Theta_0\|^2\leq C(1+t)^{5/12}.\label{3.30a}\end{equation}
So from (\ref{2.7}) and (\ref{2.9}),
\begin{eqnarray}\|\Theta-\theta_+\|^2_{L^{\infty}}\leq C\|\Theta-\theta_+\|\|\Theta_x\|^{3/4}\|\Theta_x\|^{1/4}\leq C\delta_0^{1/4}(1+t)^{-1/24}.\label{3.31a}\end{eqnarray}
 So we finish this lemma.$\Box$

We can obtain from $|(V-v_+,U,\Theta-\theta_+)|^2(x,t)\leq
C\|(V-v_+,U,\Theta-\theta_+)\|\|(V_x,U_x,\Theta_x)\|$, (\ref{2.15a})
and Lemma \ref{lem2.3} that for $x\in\R_+$,
\begin{equation}\label{3.34a}\lim_{t\to\infty}|(V,U,\Theta)|(x,t)=(v_+,0,\theta_+).\end{equation}

If
$$\|(v-V,u-U,\theta-\Theta)\|_{L^{\infty}(\R_+)}\to 0,\ \ t\to\infty,$$
we can get that the asymptotic stability results to $(v,u,\theta)$
is $(v_+,0,\theta_+)$. This stability result can be obtained at the end of the paper.

We shall prove Theorem \ref{thm1.1} by combining the local existence
and the global-in-time a priori estimates. Since the local existence
of the solution  is well known (see, for example,\cite{HMS}), we
omit it here for brevity. to prove the global existence part of
Theorem \ref{thm1.1}, it is sufficient to establish the following a
priori estimates.

\begin{pro}\label{pro2.2}{\rm(A priori estimate)} Let $(\varphi,\psi,\zeta)\in
X([0,t])$ be a solution of problem (\ref{1.12}) for some $t>0$.Set $C$ is a positive constant only depends on $C_v,\ R,\
\mu,\ \theta_{\pm},\ v_{\pm}$ and
$\|(\varphi_0,\psi_0,\zeta_0)\|_1$,
$C_0>2\left(C\|(\varphi_0,\psi_0,\zeta_0)\|_1^2+C+1\right)^{1/2}$.
If $\|(\varphi_0,\psi_0,\zeta_0)\|$ is a suitably small constant,
$$\bar{N}_1(t)=max\{m_\rho^{-1},M_{\rho},m_{\theta}^{-1},M_{\theta},\|(\varphi,\psi,\zeta)\|_1\}\leq C_0,$$
with $0<m_{\rho}=v^{-1}(x,t)\leq \rho(x,t)\leq M_{\rho}$,
$0<m_{\theta}\leq \theta(x,\tau)\leq M_{\theta}$, then
$(\varphi,\psi,\zeta)$ satisfies the a priori estimate
\begin{eqnarray}\label{3.4}
\|(\varphi,\psi,\zeta)\|_1^2+\int_0^t\left\{\|\varphi_x\|^2+\|(\psi_x,\zeta_x)\|_1^2\right\}d\tau\leq
C\|(\varphi_0,\psi_0,\zeta_0)\|_1^2+C< C_0^2/4 , \end{eqnarray} and
$\bar{N}_1(t)\leq C_0/2$.
\end{pro}

\section{Proof of Theorem \ref{thm1.1}}

   Under the preparations in last section, the main task here is to finish Proposition \ref{pro2.2} by the following lemmas.

\begin{lem}\label{lem3.2} If $C(\delta_0)>0$ is a small constant about
$\delta_0$
\begin{eqnarray*}&&\int_0^t\int_{\R_+}\Theta_x^2(\varphi^2+\zeta^2)dxd\tau\leq
C(\delta_0)\int_0^t\|(\varphi_x,\zeta_x)\|^2d\tau.
\end{eqnarray*}
\end{lem}
\pf\begin{eqnarray*}
&&\int_0^t\int_{\R_+}\Theta_x^2(\zeta^2+\varphi^2)dxd\tau\\
&&\leq\int_0^t\int_{\R_+}\Theta_x^2(\|\zeta\|\|\zeta_x\|+\|\varphi\|\|\varphi_x\|)dxd\tau\\
&&\leq C\int_0^t
(\|\zeta_x\|+\|\varphi_x\|)^2\|\Theta_x\|^{1/4}d\tau+C\int_0^t\|\Theta_x\|^{15/4}d\tau.
\end{eqnarray*}
From (\ref{2.9}) and (\ref{2.7}) we can get \begin{eqnarray*}
\int_0^t\int_{\R_+}\Theta_x^2(\zeta^2+\varphi^2)dxd\tau\leq
C(\delta_0)\int_0^t (\|\zeta_x\|+\|\varphi_x\|)^2d\tau+C(\delta_0).
\end{eqnarray*} That we finish this lemma. $\Box$

\begin{lem}\label{lem3.1}
If  $C(\delta_0)>0$ is small constant about $\delta_0$,  we can get
\begin{eqnarray}\label{4.2b}
&&\int_{\R_+}\left(\varphi^2+\psi^2+\zeta^2\right)dx+\int_0^t\
\left\|\left(\psi_x,\zeta_x\right)\right\|^2\
d\tau\nonumber\\
&&\leq C(\delta_0)+C(\delta_0)\int_0^t\
\left(\|\varphi_x\|^2+\|\psi_{xx}\|^2\right)\
d\tau+C\|(\varphi_0,\psi_0,\zeta_0)\|^2.
\end{eqnarray}
\end{lem}
\pf Set
$$\Phi(z)=z-\ln z-1,$$
$$\Psi(z)=z^{-1}+\ln z-1,$$ where $\Phi'(1)=\Phi(1)=0$ is a strictly
convex function around $z=1$. Similar to the proof in \cite{HMS},
 we deduce from (\ref{1.12}) that
\begin{eqnarray}\label{3.2}
&&\left(\frac{\psi^2}{2}+R\Theta\Phi\left(\frac{v}{V}\right)+C_v\Theta\Phi\left(\frac{\theta}{\Theta}\right)\right)_t\nonumber\\
&&\quad+\mu\frac{\Theta\psi_x^2}{v\theta}+\kappa\frac{\Theta\zeta_x^2}{v\theta^2}+H_x+Q
=\mu\left(\frac{\psi\psi_x}{v}\right)_x-F\psi-\frac{\zeta
G}{\theta},
\end{eqnarray}
where
$$H=R\frac{\zeta\psi}{v}-R\frac{\Theta\varphi\psi}{vV}+\mu\frac{U_x\varphi\psi}{vV}-\kappa\frac{\zeta\zeta_x}{v\theta}+\kappa\frac{\Theta_x\varphi\zeta}{v\theta V},$$
and
\begin{eqnarray*}
Q&=&p_+\Phi\left(\frac{V}{v}\right)U_x+\frac{p_+}{\gamma-1}\Phi\left(\frac{\Theta}{\theta}\right)U_x-\frac{\zeta}{\theta}(p_+-p)U_x-\mu\frac{U_x\varphi\psi_x}{vV}\\
&&\quad-\kappa
\frac{\Theta_x}{v\theta^2}\zeta\zeta_x-\kappa\frac{\Theta\Theta_x}{v\theta^2V}\varphi\zeta_x-2\mu\frac{U_x}{v\theta}\zeta\psi_x
+\kappa\frac{\Theta_x^2}{v\theta^2V}\varphi\zeta+\mu\frac{U_x^2}{v\theta
V}\varphi\zeta\\
&=:&\sum_{i=1}^9Q_i.
\end{eqnarray*}
Note that $p=R\theta/v$, $p_+=R\Theta/V$ and (\ref{1.7}), use
integrate by part and Cauchy-Schwarz inequality can get

\begin{eqnarray}\label{3.3}
Q_1+Q_2&=&Ra\left(\Phi\left(\frac{V}{v}\right)(\ln\Theta)_x\right)_x+\frac{Ra}{\gamma-1}\left(\Phi\left(\frac{\Theta}{\theta}\right)(\ln\Theta)_x\right)_x\nonumber\\
&&\quad-aR(\ln\Theta)_x\left(\frac{V\varphi_x\varphi-V_x\varphi^2}{Vv^2}\right)\nonumber\\
&&\quad-a\frac{p_+}{\gamma-1}(\ln\Theta)_x\left(\frac{\Theta\zeta_x\zeta-\Theta_x\zeta^2}{\Theta\theta^2}\right)\nonumber\\
&&\geq\left(p_+\Phi\left(\frac{V}{v}\right)U+\frac{p_+}{\gamma-1}\Phi\left(\frac{\Theta}{\theta}\right)U\right)_x\nonumber\\
&&\quad-\epsilon(\zeta_x^2+\varphi_x^2)-C\epsilon^{-1}\Theta_x^2(\zeta^2+\varphi^2).
\end{eqnarray}
Similarly, using $p-p_+=\frac{R\zeta-p_+\varphi}{v}$, we can get
\begin{equation}\label{3.4}
Q_3\geq\frac{R\zeta-p_+\varphi}{v}(\frac{\zeta}{\theta}U_x)\geq\left(\frac{R\zeta^2U}{v\theta}-\frac{p_+\zeta\varphi
U}{\theta
v}\right)_x-C(\delta_0)(\zeta_x^2+\varphi_x^2)-C^{-1/2}(\delta_0)\Theta_x^2(\zeta^2+\varphi^2).
\end{equation}
And
\begin{eqnarray}\label{3.5}
(Q_4+Q_7)+(Q_5+Q_6+Q_8)+Q_9&\geq&
-CC^{-1/2}(\delta_0)(\ln\Theta)_{xx}^2-C^{1/2}(\delta_0)\psi_x^2\nonumber\\
&&-C^{1/2}(\delta_0)\zeta_x^2-CC^{-1/2}(\delta_0)\Theta_x^2(\zeta^2+\varphi^2)\nonumber\\
&&-CC^{-1/2}(\delta_0)|(\ln\Theta)_{xx}|^2(\zeta^2+\varphi^2).
\end{eqnarray}
At the end we use the definition of $F$ and $G$ in (\ref{1.10}) then
combine with the general inequality skills as above to get
\begin{eqnarray}\label{3.6}
-F\psi-G\frac{\zeta}{\theta}&=&-\frac{\kappa a(\gamma-1)-\mu
p_+\gamma}{R\gamma}\left(\frac{(\ln\Theta)_{xx}}{\Theta}\right)_x\psi\nonumber\\
&&\quad+\frac{\mu p_+}{R\Theta}\left(\frac{\kappa
(\gamma-1)}{R\gamma}(\ln\Theta)_{xx}\right)^2\frac{\zeta}{\theta}\nonumber\\
&\leq&-\frac{\kappa a(\gamma-1)-\mu
p_+\gamma}{R\gamma}\left(\frac{(\ln\Theta)_{xx}}{\Theta}\psi\right)_x+\frac{\kappa
a(\gamma-1)-\mu
p_+\gamma}{R\gamma}\frac{(\ln\Theta)_{xx}}{\Theta}\psi_x\nonumber\\
&&\quad+\frac{\mu p_+}{R\Theta}\left(\frac{\kappa
(\gamma-1)}{R\gamma}(\ln\Theta)_{xx}\right)^2\frac{\zeta}{\theta}\nonumber\\
&\leq&-\frac{\kappa a(\gamma-1)-\mu
p_+\gamma}{R\gamma}\left(\frac{(\ln\Theta)_{xx}}{\Theta}\psi\right)_x+C^{1/2}(\delta_0)\psi_x^2+CC^{-1/2}(\delta_0)(\ln\Theta)_{xx}^2.
\end{eqnarray}
Integrating (\ref{3.3}) to (\ref{3.6}) over $\R\times(0,t)$ , using
Lemma \ref{lem2.3} and the boundary condition about
$(\varphi,\psi,\zeta)$ of (\ref{1.12}) to estimate the terms
$\mu\left(\frac{\psi\psi_x}{v}\right)_x$,
$\left(\frac{(\ln\Theta)_{xx}\psi}{\Theta}\right)_x$ and $H_x$, in
the end combine with Cauchy-Schwarz inequality  we know   that for a
small $C(\delta_0)>0$ which is about $\delta_0$, we have
\begin{eqnarray}\label{3.7}
&&\int_{\R_+}\left(R\Theta\Phi\left(\frac{v}{V}\right)+\frac{1}{2}\psi^2+C_v\Theta\Phi\left(\frac{\theta}{\Theta}\right)\right)dx+\int_0^t\
\left\|\left(\psi_x/(\sqrt{v\theta}),\zeta_x/(\theta\sqrt{v})\right)\right\|^2\
d\tau\nonumber\\
&&\leq
C^{-1/2}(\delta_0)\left\{\int_0^t\int_0^{+\infty}\Theta_x^2(\varphi^2+\zeta^2)\
dxd\tau+\|\Theta_{0x}\|^2\right\}+C\left\{C^{1/2}(\delta_0)\int_0^t\ \|\varphi_x\|^2\ d\tau+\|(\varphi_0,\psi_0,\zeta_0)\|^2\right\}\nonumber\\
&&\qquad+C^{-1/2}(\delta_0)\int_0^t\psi^2(0,\tau)d\tau+C^{1/2}(\delta_0)\int_0^t\psi_x^2(0,\tau)
d\tau+C^{-1/2}(\delta_0)\int_0^t(\ln\Theta)_{xx}^2(0,\tau)d\tau+C(\delta_0).
\end{eqnarray}

Using the definition about $\psi(0,t)$ in (\ref{1.12}), then combine
with (\ref{1.9})$_5$,Cauchy-Schwarz inequality and Lemma
\ref{lem2.3} we can get
\begin{eqnarray}
&&\int_0^t\psi^2(0,\tau)d\tau+\int_0^t(\ln\Theta)_{xx}^2(0,\tau)d\tau\nonumber\\
&&\leq C\int_0^t(\|(\ln\Theta)_x\|^{1/2}\|(\ln\Theta)_x\|^{7/2}+\|(\ln\Theta)_{xx}\|^{1/15}\|(\ln\Theta)_{xx}\|^{19/15})d\tau\nonumber\\
&&\quad+\int_0^t(\|(\ln\Theta)_{xx}\|^2+\|\partial_x^3(\ln\Theta)\|^2)d\tau\nonumber\\
&&\leq C(\delta_0).\label{3.8}
\end{eqnarray}
Because $$\int_0^t\psi^2_x(0,\tau)d\tau\leq C\int_0^t(\|\psi_x\|^2+\|\psi_{xx}\|^2)d\tau,$$
combine with (\ref{3.8}) and Lemma
\ref{lem3.2},   (\ref{3.7}) can be change to (\ref{4.2b}).$\Box$

\begin{lem}\label{lem3.3} If  $\epsilon$ is a positive constant, $\delta>0$ stands for a small constant about
$\|(\varphi_0,\psi_0,\zeta_0)\|$ and $\delta_0$, we can get
\begin{eqnarray*}
&&\|(\varphi,\psi,\zeta)\|^2+\|(\psi_x,\zeta_x)\|^2+\int_0^t\|(\psi_{xx},\zeta_{xx})\|^2d\tau\\
&&\leq
C\|(\psi_{0x},\zeta_{0x})\|^2+C\delta+C\epsilon^{-1}\int_0^t\|\varphi_x\|^2d\tau.
\end{eqnarray*}.
\end{lem}

\pf First to get the estimate of $\|\psi_x(t)\|$ ,multiply both side
of  (\ref{1.12})$_2$ by $\psi_{xx}$ to get
\begin{eqnarray*}
&&\left(\frac{\psi_x^2}{2}\right)_t+\mu\frac{\psi_{xx}^2}{v}
=\mu\frac{\psi_x v_x}{v^2}\psi_{xx}+\mu\left(\frac{U_x\varphi}{v
V}\right)_x\psi_{xx}\\
&&\quad-R\left(\frac{\Theta\varphi}{v
V}\right)_x\psi_{xx}+R\left(\frac{\zeta}{v}\right)_x\psi_{xx}+F\psi_{xx}+(\psi_t\psi_x)_x:=\sum_{i=1}^6I_i.
\end{eqnarray*}
When we integrate it in $\R_+\times (0,t)$  , we
get
\begin{eqnarray}
&&\|\psi_x(t)\|^2+\int_0^t\|\psi_{xx}(\tau)\|^2d\tau\nonumber\\
&&\quad\leq
C\|\psi_{0x}\|^2+C\sum_{i=1}^6\left|\int_0^t\int_0^\infty
I_idxd\tau\right|.\label{3.9}
\end{eqnarray}

Now we deal with   $\iint|I_i|dxd\tau$ in the right side of
(\ref{3.9}). Using $v=\varphi+V$ , $R\Theta/V=p_+$ , Lemma
\ref{lem2.3}, we can get
\begin{eqnarray}
&&\int_0^t \int_0^{+\infty}|I_1| dx d\tau\leq C\int_0^t
\int_0^{+\infty}|V_x|
|\psi_{x}||\psi_{xx}|dxd\tau+C\int_0^t\int_0^{+\infty}
|\varphi_x||\psi_{x}||\psi_{xx}|dx d\tau\nonumber\\
&&\quad\leq C\int_0^t\|V_x\|\|\psi_x\|_{L^\infty}\|\psi_{xx}\|
d\tau+C\int_0^t \|\psi_x\|_{L^\infty}\|\varphi_x\|\|\psi_{xx}\|
d\tau\nonumber\\
&&\quad\leq C^{1/2}(\delta_0)\int_0^t\|\psi_{xx}\|^2d\tau+
C^{-1/2}(\delta_0)\int_0^t\|\psi_x\|^2 \|V_x\|^4d\tau+C\int_0^t
\|\psi_x\|^{1/2}\|\varphi_x\|\|\psi_{xx}\|^{3/2}
d\tau\nonumber\\
&&\quad\leq C(\delta_0)\int_0^t\|\psi_{xx}\|^2d\tau+
C(\delta_0)\int_0^t\|\psi_x\|^2
d\tau+C^{-1/4}(\delta_0)\sup_{t}\|\varphi_x\|^4\int_0^t
\|\psi_x\|^2 d\tau .\label{3.10}
\end{eqnarray}
Because $N_1(t)\leq C_0$ and $\|(\varphi_0,\psi_0,\zeta_0)\|$ is small, we can get from
 Lemma \ref{lem3.1} that
$$\int_0^t\int_{\R_+}|I_1|dxd\tau\leq C(\delta_0)\int_0^t\|\psi_{xx}\|^2d\tau
+C\delta\int_0^t\|\varphi_x\|^2d\tau+C\delta\int_0^t\|\psi_{xx}\|^2d\tau+\delta.$$
Next we use the definition of $(V,U,\Theta)$ (see (\ref{1.7}),
(\ref{1.9}) and (\ref{2.15a})), Cauchy-Schwarz inequality and Lemma \ref{lem2.3} to get
\begin{eqnarray}
&&\int_0^t\int_0^\infty|I_2|dxd\tau\nonumber\\
&&\quad\leq
C\int_0^t\int_0^\infty\left(|U_{xx}||\varphi|+|U_x||\varphi_x|+|U_x||V_x||\varphi|+|U_x||\varphi||\varphi_x|\right)|\psi_{xx}|dxd\tau\nonumber\\
&&\quad\leq
C^{1/2}(\delta_0)\int_0^t\|\psi_{xx}\|^2d\tau+\frac{C}{C^{1/2}(\delta_0)}\int_0^t\|\varphi\|_{L^\infty}^2\|U_{xx}\|^2d\tau
+\frac{C}{C^{1/2}(\delta_0)}\int_0^t\|U_x\|_{L^\infty}^2\|\varphi_x\|^2d\tau\nonumber\\
&&\qquad+\frac{C}{C^{1/2}(\delta_0)}\int_0^t\|\varphi\|_{L^\infty}^2\|V_x\|^2\|U_x\|_{L^\infty}^2d\tau
+\frac{C}{C^{1/2}(\delta_0)}\int_0^t\|\varphi\|_{L^\infty}^2\|U_x\|_{L^\infty}^2\|\varphi_x\|^2d\tau\nonumber\\
&&\quad\leq
C^{1/2}(\delta_0)\int_0^t\|\psi_{xx}\|^2d\tau+C(\delta_0)+C(\delta_0)\int_0^t\|\varphi_x\|^2d\tau.\label{3.11}
\end{eqnarray}
The same as (\ref{3.10}) and (\ref{3.11}), we use Lemma
\ref{lem3.2}, the definition of $F$ in (\ref{1.10}),  Lemma
\ref{lem2.3} and (\ref{2.15a}) we can get the estimates about $I_3$ to
$I_5$ as following.
\begin{eqnarray}
&&\int_0^t\int_0^\infty(|I_3|+|I_4|+|I_5|)dxd\tau\nonumber\\
&&\leq
C\int_0^t\int_0^\infty\left(|\Theta_{x}||\varphi|+|\Theta||\varphi_x|+|\Theta||V_x||\varphi|
+|\Theta||\varphi||\varphi_x|\right)|\psi_{xx}|dxd\tau\nonumber\\
&&\quad+C\int_0^t\int_0^\infty\left(|\zeta_x|+|\zeta||V_x|+|\zeta||\varphi_x|\right)|\psi_{xx}|dxd\tau\nonumber\\
&&\quad+C^{1/2}(\delta_0)\int_0^t\|\psi_{xx}\|^2d\tau+\frac{C}{C^{1/2}(\delta_0)}\int_0^t\|F\|^2d\tau\nonumber\\
&&\leq
C^{1/2}(\delta_0)\int_0^t\|\psi_{xx}\|^2d\tau+CC^{-1/2}(\delta_0)\int_0^t\|\varphi_x\|^2d\tau
+\frac{C}{C^{1/2}(\delta_0)}\int_0^t\|\varphi_x\|^2d\tau\nonumber\\
&&\quad+\frac{C}{C^{1/2}(\delta_0)}\int_0^t\int_0^\infty
V_x^2\varphi^2dx
d\tau+C^{1/2}(\delta_0)\int_0^t\|\psi_{xx}\|^2d\tau+\frac{C}{C^{1/2}(\delta_0)}\int_0^t\int_0^\infty\left(\zeta_x^2+V_x^2\zeta^2\right)dxd\tau\nonumber\\&&\quad+\frac{C}{C^{1/2}(\delta_0)}\sup_t\|(\varphi,\zeta)\|\|(\varphi_x,\zeta_x)\|\int_0^t\|\varphi_x\|^2d\tau
+C^{1/2}(\delta_0)\int_0^t\|\psi_{xx}\|^2d\tau+C(\delta_0).\label{3.12}
\end{eqnarray}
Because $N_1(t)\leq C_0$ and $\|(\varphi_0,\psi_0,\zeta_0)\|$ is small, we can get from
 Lemma \ref{lem3.1} that
$$\frac{C}{C^{1/2}(\delta_0)}\sup_t\|(\varphi,\zeta)\|\|(\varphi_x,\zeta_x)\|\int_0^t\|\varphi_x\|^2d\tau\leq C\delta\int_0^t\|\varphi_x\|^2d\tau+C\delta\int_0^t\|\psi_{xx}\|^2d\tau+\delta,$$

 Therefore
\begin{eqnarray*}
&&\int_0^t\int_0^\infty(|I_3|+|I_4|+|I_5|)dxd\tau\\
&&\leq
\epsilon\int_0^t\|\psi_{xx}\|^2d\tau+\frac{C}{\epsilon}\int_0^t\int_0^\infty\
(\zeta_x^2+\varphi_x^2)dxd\tau
+C\delta\int_0^t\|\psi_{xx}\|^2d\tau+C(\delta_0)+\delta.\end{eqnarray*}

At last we use integration by parts to the term about $I_6$ to get
\begin{eqnarray}
&&\left|\int_0^t\int_0^\infty
I_6dxd\tau\right|=\left|\int_0^t(\psi_t\psi_x)(0,\tau)d\tau\right|\leq
\frac{C}{C^{1/2}(\delta_0)}\int_0^t\psi_x^2(0,\tau)d\tau+C^{1/2}(\delta_0)\int_0^t\psi_\tau^2(0,\tau)d\tau\nonumber\\
&&\leq C^{-1/2}(\delta_0)\int_0^t \|\psi_x\|^2\ d\tau+1/16\int_0^t
\|\psi_{xx}\|^2\
d\tau+C^{1/2}(\delta_0)\int_0^t\psi_\tau^2(0,\tau)d\tau.\label{3.13}
\end{eqnarray}
Using the definition of $U$ in (\ref{1.7}), $\psi=u-U$ and
(\ref{2.29}) to get
\begin{eqnarray}
\psi_t(0,t)&=&-\frac{k(\gamma-1)}{\gamma
R}(\ln\Theta)_{xt}(0,t)\nonumber\\
&&=-a\frac{k(\gamma-1)}{\gamma
R}\partial_x\left(\frac{(\ln\Theta)_{xx}}{\Theta}\right)(0,t).
\label{3.14}
\end{eqnarray}

Combine with Lemma \ref{lem2.3} and (\ref{2.32b}) we get
\begin{eqnarray}
\int_0^t\|(\ln\Theta)_{x\tau}\|_{L^\infty}^2(0,\tau)d\tau&\leq&
C\int_0^t\|(\ln\Theta)_{xxx}\|^2d\tau+C\int_0^t\|\partial_x^4(\ln\Theta)\|^2d\tau\nonumber\\
&&\quad+C\int_0^t\|(\ln\Theta)_x\|_{L^{\infty}}\|(\ln\Theta)_{xx}\|_{L^{\infty}}\leq C
.\label{3.15}
\end{eqnarray}
So combine with Lemma \ref{lem3.1}, (\ref{3.13}), (\ref{3.14}) and (\ref{3.15}) we get
\begin{eqnarray}
&&\left|\int_0^t\int_0^\infty I_6dxd\tau\right|\nonumber\\
&&\leq
 \frac{C}{C^{1/2}(\delta_0)}\int_0^t \|\psi_x\|^2\
d\tau+1/16\int_0^t \|\psi_{xx}\|^2\ d\tau+C(\delta_0)\nonumber\\
&&\leq \delta+C\delta\int_0^t\|\varphi_x\|^2+1/16\int_0^t \|\psi_{xx}\|^2\ d\tau.\label{3.16}
\end{eqnarray}

In all, there exist a small constant $\delta$ which is about
$\|(\varphi_0,\psi_0,\zeta_0)\|$ and $\delta_0$, such that
\begin{eqnarray}\label{3.17}
&&\int_0^t\ \int_0^{+\infty}\ \sum_{i=1}^6|I_i|\ dx\ d\tau\nonumber\\
&&\leq \int_0^t\ 1/2\|\psi_{xx}\|^2d\tau+C\epsilon^{-1}\int_0^t\
\|\varphi_x\|^2\ d\tau+C\delta.
\end{eqnarray}

So (\ref{3.9}) can be change to
\begin{eqnarray}
&&\|\psi_x(t)\|^2+\int_0^t\|\psi_{xx}(\tau)\|^2d\tau\nonumber\\
&&\leq
C\delta+C\epsilon^{-1}\int_0^t\|\varphi_x\|^2d\tau+C\|\psi_{0x}\|^2.
\label{3.18}
\end{eqnarray}

The estimate about  $\|\zeta_x\|$ is similar to $\|\psi_x\|$, use
(\ref{1.12})$_3$ multiply $\zeta_{xx}$ then integrate in
$Q_t=\R_+\times(0,t)$ to get
\begin{eqnarray}
&&\|\zeta_x\|^2+\int_0^t
\|\zeta_{xx}\|^2\ d\tau\nonumber\\
&&\leq C \|\zeta_{0x}\|^2+
C\int_0^t\int_0^\infty\left(\psi_x^2+\zeta^2\psi_x^2+\zeta^2U_x^2+U_x^2\varphi^2\right)dxd\tau\nonumber\\
&&
\quad+C\int_0^t\int_0^{+\infty}|\zeta_x|(|\varphi_x|+|V_x|)|\zeta_{xx}|dxd\tau+C\int_0^t\int_0^\infty\left|\left(\frac{\Theta_x\varphi}{vV}\right)_x\right|^2dxd\tau
\nonumber\\
&&\quad+C\int_0^t\int_0^{+\infty}(U_x^4+\psi_x^4)dxd\tau+C\int_0^t\|G\|^2d\tau\nonumber\\
&&=:C\|\zeta_{0x}\|^2+\sum_{i=1}^5J_i.\label{3.19}
\end{eqnarray}
Use the same method as (\ref{3.10})--(\ref{3.13}) and combine with
Lemma \ref{lem3.1}
$$
J_1\leq
C(1+N^2(t))\int_0^t\|\psi_x\|^2d\tau+CN^2(t)\int_0^t\|U_x\|^2d\tau
\leq
C\delta\left(1+\int_0^t\|\psi_{xx}\|^2d\tau+\int_0^t\|\varphi_x\|^2d\tau\right)+C(\delta_0).
$$
Again use the same method as (\ref{3.10})--(\ref{3.13}) and combine
with Lemma \ref{lem3.1}
\begin{eqnarray*}
J_2&\leq& C\int_0^t\|\zeta_x\|_{L^\infty}\|\varphi_x\|\|\zeta_{xx}\|d\tau+C\int_0^t\|V_x\|\|\zeta_x\|_{L^\infty}\|\zeta_{xx}\|d\tau\\
&\leq&C\int_0^t\|\zeta_x\|^{1/2}\|\zeta_{xx}\|^{3/2}\|\varphi_x\|d\tau
+1/16\int_0^t\|\zeta_{xx}\|^2d\tau+C(\delta_0)\int_0^t\|\zeta_x\|^2d\tau\\
&\leq&1/8\int_0^t\|\zeta_{xx}\|^2d\tau+C(\delta_0)\int_0^t\|\zeta_x\|^2d\tau
+C\sup_t\|\varphi_x\|^4\int_0^t\|\zeta_x\|^2d\tau\\
&\leq&C\delta\left(1+\int_0^t\|\psi_{xx}\|^2d\tau+\int_0^t\|\varphi_x\|^2d\tau\right)+C(\delta_0)+1/8\int_0^t\|\zeta_{xx}\|^2d\tau.
\end{eqnarray*}

Because
\begin{eqnarray*}&&\left|\left(\frac{\Theta_x\varphi}{vV}\right)_x\right|^2\\
&&=|\frac{\Theta_{xx}\varphi}{vV}+\frac{\Theta_x\varphi_x}{vV}+\frac{\Theta_x\varphi}{vV}(-\frac{V_x+\varphi_x}{v^2}-\frac{V_x}{V^2})|^2\\
&&\leq
C\Theta_{xx}^2\varphi^2+C\Theta_x^2\varphi_x^2+C\Theta_x^2V_x^2\varphi^2+C\Theta_x^2\varphi^2\varphi_x^2,
\end{eqnarray*}
combine with $R\Theta/V=p_+$, use the same method as
(\ref{3.10})--(\ref{3.13}) to get
\begin{eqnarray*}
J_3&\leq& C\int_0^t\|\varphi\|_{L^\infty}^2\|\Theta_{xx}\|^2d\tau
+C\int_0^t\|\Theta_x\|_{L^\infty}^2\|\varphi_x\|^2d\tau\\
&&\quad+C\int_0^t\ \int_0^{+\infty}\ \Theta_x^2V_x^2\varphi^2\ dx\ d\tau\nonumber\\
&\leq&C(\delta_0)\int_0^t\|\varphi_{x}\|^2d\tau+C(\delta_0) .
\end{eqnarray*}
Use the definition $U$ and similar as (\ref{3.10}) (\ref{3.11})
that we combine with Lemma \ref{lem2.3} to get
\begin{eqnarray*} J_4&\leq&
C(\delta_0)+C\int_0^t\|\psi_{x}\|_{L^\infty}^2\|\psi_x\|^2d\tau\\
&\leq& C(\delta_0)+C\int_0^t\|\psi_x\|^3\|\psi_{xx}\|d\tau\\
&\leq&
C(\delta_0)+C\int_0^t\left(\|\psi_x\|^2\|\psi_x\|^4+1/16\|\psi_{xx}\|^2\right)d\tau\\
&\leq&C\delta
N^4(t)\left(1+\int_0^t\|\psi_{xx}\|^2d\tau+\int_0^t\|\varphi_x\|^2d\tau\right)+C(\delta_0)+1/16\int_0^t\|\psi_{xx}\|^2d\tau.
\end{eqnarray*}
Use the definition $G$ in (\ref{2.2}) combine with Lemma
\ref{lem2.3}
\begin{eqnarray*}
J_5=C\int_0^t\|G\|^2d\tau\leq C(\delta_0).
\end{eqnarray*}
Use the results from $J_1$ to $J_5$, the inequality (\ref{3.19}) can
be change to
\begin{eqnarray}
&&\|\zeta_x\|^2+\int_0^t\ \|\zeta_{xx}\|^2\ d\tau\nonumber\\
&&\leq
C\|\zeta_{0x}\|^2+C\delta\left(1+\int_0^t\|\psi_{xx}\|^2d\tau+\int_0^t\|\varphi_x\|^2d\tau\right)+1/16\int_0^t\|\psi_{xx}\|^2d\tau.\label{3.20}
\end{eqnarray}

 In fact when combine  with Lemma \ref{lem3.1}--\ref{lem3.2}, (\ref{3.18}) and (\ref{3.20}), it is easy to get
\begin{eqnarray*}
&&\|(\varphi,\psi,\zeta)\|^2+\|(\psi_x,\zeta_x)\|^2+\int_0^t\|(\psi_{xx},\zeta_{xx})\|^2d\tau\\
&&\leq
C\|(\psi_{0x},\zeta_{0x})\|^2+C\delta+C\epsilon^{-1}\int_0^t\|\varphi_x\|^2d\tau.
\end{eqnarray*}
$\Box$

\begin{lem}\label{lem3.4} For a small $C(\delta)>0$  stands for constant about $\|(\varphi_0,\psi_0,\zeta_0)\|$ and $\delta_0$, and $C(\delta_0)>0$ is a small constant about
$\delta_0$  , we can get
\begin{eqnarray}
\|\varphi_x\|^2+\int_0^t\|\varphi_x\|^2 d\tau&\leq&
C\|(\varphi_{0x},\psi_{0x},\zeta_{0x})\|^2+C\delta .\label{3.21}
\end{eqnarray}
\end{lem}
\pf Set $\bar{v}=\frac{v}{V}$, take it into (\ref{1.12})$_1$,
(\ref{1.12})$_2$  ($p=R \theta/v$) to get
$$
\psi_t+p_x=\mu\left (\frac{\bar{v}_x}{\bar{v}}\right)_t-F.
$$
Both sides of last equation multiply $\bar v_x/\bar v$ to get
\begin{eqnarray}
&&\left(\frac{\mu}{2}\left (\frac{\bar v_x}{\bar
v}\right)^2-\psi\frac{\bar v_x}{\bar v}\right)_t
+\frac{R\theta}{v}\left (\frac{\bar v_x}{\bar
v}\right)^2+\left(\psi\frac{\bar v_t}{\bar v}\right)_x\nonumber\\
&&\quad=\frac{\psi_x^2}{v}+U_x\left
(\frac{1}{v}-\frac{1}{V}\right)\psi_x+\frac{R\zeta_x}{v}\frac{\bar
v_x}{\bar v}-\frac{R\theta}{v}\left
(\frac{1}{\Theta}-\frac{1}{\theta}\right)\Theta_x\frac{\bar
v_x}{\bar v}+F\frac{\bar v_x}{\bar v}.\label{3.22}
\end{eqnarray}
Because $v|_{x=0}=V|_{x=0}=v_-$, we can get
$$
\left(\psi\frac{\bar v_t}{\bar v}\right)\Big|_{x=0}=0.$$ On the
other hand if we integrate (\ref{3.22}) in $R_+\times(0,t)$, combine
with Cauchy-Schwartz inequality,
 (\ref{3.22}) is changed to
\begin{eqnarray*}
&&\int_{\R_+}\left(\frac{\mu}{2}\left (\frac{\bar v_x}{\bar
v}\right)^2-\psi\frac{\bar v_x}{\bar
v}\right)dx+\int_0^t\int_{\R_+}\frac{R\theta}{v}\left (\frac{\bar
v_x}{\bar v}\right)^2dxd\tau\nonumber\\
&&\leq C^{-1/2}(\delta_0)\left(\int_0^t\ \|(\zeta_x,\psi_x)\|^2\
d\tau+\int_0^t\ \int_0^{+\infty}\ \Theta_x^2(\varphi^2+\zeta^2)\ dx\
d\tau\right)\\&&+C\int_0^t\int_0^{+\infty}U_x^2\varphi^2
dxd\tau+C\int_0^t\int_0^{+\infty}|F|^2dxd\tau+1/2\int_0^t\|\frac{\sqrt{R\theta}}{\sqrt{v}}\frac{\bar{v}_x}{\bar{v}}\|^2\
d\tau.\end{eqnarray*} Furthermore, because
$C_1(\varphi_x^2)-C_2V_x^2\leq (\frac{\bar{v}_x}{\bar{v}})^2\leq
C_3\varphi_x^2+C_4V_x^2$ ($C_1,C_2,C_3,C_4$ stands for constants
about $v$), combine with Lemma \ref{lem3.2}, Lemma \ref{lem3.1} and
Lemma \ref{lem3.3} we can get
\begin{eqnarray}\label{3.25}\int_0^t\|\varphi_x\|^2\
d\tau+\|\varphi_x\|^2&\leq&
C\|(\varphi_{0x},\psi_{0x},\zeta_{0x})\|^2+C\delta.
\end{eqnarray}
 So we finish this lemma.$\Box$

 From Lemma \ref{lem3.1} to Lemma \ref{lem3.4} we know when $\delta_0$ and $\|(\varphi_0,\psi_0,\zeta_0)\|$ suitably small there exist a suitably
 small positive
 constant $\delta$
such that
$$\|(\varphi,\psi,\zeta)\|^2+\int_0^t\|(\psi_x,\zeta_x)\|^2d\tau\leq
 C\delta,$$
 Then we can get $$|v-V|^2\leq \|\varphi\|\|\varphi_x\|\leq C\delta,$$
 which means $C_5\leq |v|\leq
 C_6$. Use this result to Lemma \ref{lem3.3}, we can get there exist a positive constant $C$ independent of $v(x,t)$,$u(x,t)$ and $\theta(x,t)$ such that $$\|(\varphi_x,\psi_x,\zeta_x)\|^2+\int_0^t\|(\psi_{xx},\zeta_{xx})\|^2\leq
 C\|(\varphi,\psi,\zeta)\|_1^2.$$ Similar as the estimates for the upper and lower of $v$, when we combine $\|\zeta_x\|\leq C$ with Lemma \ref{lem3.1}, we can obtain $C_7\leq |\theta|\leq C_8$ . Here $C_5$, $C_6,$ $C_7$ and
 $C_8$ are constants independent of $v(x,t)$, $u(x,t)$ and $\theta(x,t)$ .So we finish
 Proposition 2.2 .

To finish Theorem \ref{thm1.1} now we will proof $\sup_{x\in
\R_+}|(\varphi,\psi,\zeta)|\to 0,\ as\ t\to \infty.$

Because $\int_0^{+\infty}\partial_x$(\ref{1.12})$_1\times
2\varphi_x\ dx$
 equals to
\begin{equation}
0=2\int_0^\infty\varphi_x\psi_{xx}dx-\frac{d}{dt}\|\varphi_x\|^2,\label{3.26}
\end{equation}
use  Cauchy-Schwarz inequality we get
$$2\int_0^\infty\varphi_x\psi_{xx}dx \leq C\left(\|\varphi_x\|^2+\|\psi_{xx}\|^2\right),$$
again using Lemma \ref{lem3.3}--\ref{lem3.4} and (\ref{3.26}), then
 we get
\begin{eqnarray}
&&\int_0^\infty\left|\frac{d}{dt}\|\varphi_x(t)\|^2\right|dt\nonumber\\
&&\leq C\int_0^\infty\left(\|\varphi_x\|^2+\|\psi_{xx}\|^2\right)dt\nonumber\\
&&\leq C\|(\varphi_0,\psi_0,\zeta_0)\|_1^2+C\delta.\label{3.28}
\end{eqnarray}
Similar as above, from Lemma \ref{lem3.1}$-$\ref{lem3.4} and combine
with Sobolev inequality we get
\begin{equation}
\int_0^\infty\left(\left|\frac{d}{dt}\|\psi_x(t)\|^2\right|+\left|\frac{d}{dt}\|\zeta_x(t)\|^2\right|\right)d\tau\leq
C\|(\varphi_0,\psi_0,\zeta_0)\|_1^2+C\delta.\label{3.29}
\end{equation}

It means
$$
\|(\varphi,\psi,\zeta)(t)\|_{L^\infty}^2\leq
2\|(\varphi,\psi,\zeta)(t)\|\|(\varphi_x,\psi_x,\zeta_x)(t)\|\to0\quad\mbox{when}\quad
t\to\infty.
$$
 Now when we combine with (\ref{3.34a}) we finish the theorem.

\end{document}